\documentclass[10pt]{amsart}
\usepackage{graphicx,amscd,color,amsmath,amsfonts,amssymb,geometry,soul}

\newtheorem{theorem}{Theorem}
\theoremstyle{plain}

\newtheorem{remark}{Remark}

\numberwithin{equation}{section}
\begin{document}
\title[A note on absolutely summing multilinear operators]{A note on
absolutely summing multilinear operators}
\author[Nogueira]{T. Nogueira}
\address[T. Nogueira]{Departamento de Matem\'{a}tic,a\\
\indent Universidade Federal da Para\'{\i}ba\\
\indent 58.051-900 - Jo\~{a}o Pessoa, Brazil.}
\email{tonykleverson@gmail.com}
\author[Pellegrino]{D. Pellegrino}
\address[D. Pellegrino]{Departamento de Matem\'{a}tic,a\\
\indent Universidade Federal da Para\'{\i}ba\\
\indent 58.051-900 - Jo\~{a}o Pessoa, Brazil.}
\email{dmpellegrino@gmail.com and pellegrino@pq.cnpq.br}
\thanks{D. Pellegrino is supported by CNPq}
\keywords{multiple summing operators, absolutely summing operators,
Bohnenblust--Hille inequality}

\begin{abstract}
We obtain some optimal estimates for multilinear forms on $\ell _{p}$ spaces.
\end{abstract}

\maketitle


\section{Introduction}

\bigskip A result due to Zalduendo (\cite[Corollary 1]{zalduendo}) asserts
that if $m\geq 2$ is a positive integer and $p>m$, then
\begin{equation}
\left( \sum\limits_{j=1}^{n}\left\vert T\left( e_{j},...,e_{j}\right)
\right\vert ^{\frac{p}{p-m}}\right) ^{\frac{p-m}{p}}\leq \left\Vert
T\right\Vert  \label{zz}
\end{equation}%
for all $m$-linear forms $T:\ell _{p}^{n}\times \cdots \times \ell
_{p}^{n}\rightarrow \mathbb{K}$ (here and henceforth $\mathbb{K=R}$ or $%
\mathbb{C}$) and all positive integers $n$, and the exponent $\frac{p}{p-m}$
is optimal. In this note we investigate what happens if $1<p\leq m$ and what
happens if we still consider $p>m$ but with exponents $s\neq \frac{p}{p-m}.$
Our results are shown to be optimal, except for the case $1<p\leq 2.$ Our
result is the following:

\begin{theorem}
Let $m\geq 2$ be a positive integer.

(a) If $s\geq 1$ and $p>m$, then
\begin{equation}
\left( \sum\limits_{j=1}^{n}\left\vert T\left( e_{j},...,e_{j}\right)
\right\vert ^{s}\right) ^{\frac{1}{s}}\leq \left\Vert T\right\Vert n^{\max
\left\{ \frac{m}{p}+\frac{1}{s}-1,0\right\} }  \label{nt}
\end{equation}%
for all $m$-linear forms $T:\ell _{p}^{n}\times \cdots \times \ell
_{p}^{n}\rightarrow \mathbb{K}$ and all positive integers $n$, and the
exponent is optimal.

(b) If $s\geq 1$ and $2\leq p\leq m$, then
\begin{equation*}
\left( \sum\limits_{j=1}^{n}\left\vert T\left( e_{j},...,e_{j}\right)
\right\vert ^{s}\right) ^{\frac{1}{s}}\leq \left\Vert T\right\Vert
\end{equation*}%
for all $m$-linear forms $T:\ell _{p}^{n}\times \cdots \times \ell
_{p}^{n}\rightarrow \mathbb{K}$ and all positive integers $n$, and the
result is optimal.

(c) If $s\geq \frac{2}{m}$ and $1<p<2$, then
\begin{equation*}
\left( \sum\limits_{j=1}^{n}\left\vert T\left( e_{j},...,e_{j}\right)
\right\vert ^{s}\right) ^{\frac{1}{s}}\leq \left\Vert T\right\Vert n^{\frac{%
2ms+2p-spm}{2sp}}
\end{equation*}%
for all $m$-linear forms $T:\ell _{p}^{n}\times \cdots \times \ell
_{p}^{n}\rightarrow \mathbb{K}$ and all positive integers $n$.
\end{theorem}

\begin{remark}
For sums with multiple indices a similar result was recently obtained in
\cite{archiv}.
\end{remark}

\section{\textbf{Proof }}

\textbf{Proof of (a). } From now on, $T:\ell _{p}^{n}\times \dots \times
\ell _{p}^{n}\rightarrow \mathbb{K}$ is an $m$-linear form. \ If $s<\frac{p}{%
p-m}$, let $x$ be such that
\begin{equation*}
\frac{1}{s}=\frac{1}{p/(p-m)}+\frac{1}{x}.
\end{equation*}%
Using H\"{o}lder's inequality we have
\begin{align*}
& \left( \sum\limits_{j=1}^{n}\left\vert T\left( e_{j},...,e_{j}\right)
\right\vert ^{s}\right) ^{\frac{1}{s}} \\
& \leq \left( \sum_{j=1}^{n}\left\vert T(e_{j},\ldots ,e_{j})\right\vert ^{%
\frac{p}{p-m}}\right) ^{\frac{p-m}{p}}\cdot \left( \sum_{j=1}^{n}\left\vert
1\right\vert ^{x}\right) ^{\frac{1}{x}} \\
& =\left( \sum_{j=1}^{n}\left\vert T(e_{j},\ldots ,e_{j})\right\vert ^{\frac{%
p}{p-m}}\right) ^{\frac{p-m}{p}}\cdot n^{\frac{1}{x}} \\
& \leq \Vert T\Vert \cdot n^{\frac{1}{s}-\frac{p-m}{p}} \\
& =\Vert T\Vert \cdot n^{\max \left\{ \frac{m}{p}+\frac{1}{s}-1,0\right\} }
\end{align*}

\bigskip Now suppose that $s\geq \frac{p}{p-m}.$ From the inclusion theorem
for $\ell _{p}$ spaces we have, invoking Zalduendo's estimate (\ref{zz})%
\begin{eqnarray*}
\left( \sum\limits_{j=1}^{n}\left\vert T\left( e_{j},...,e_{j}\right)
\right\vert ^{s}\right) ^{\frac{1}{s}} &\leq &\left( \sum_{i_{1},\ldots
,i_{m}=1}^{n}\left\vert T(e_{j},\ldots ,e_{j})\right\vert ^{\frac{p}{p-m}%
}\right) ^{\frac{p-m}{p}} \\
&\leq &\Vert T\Vert \\
&=&\Vert T\Vert \cdot n^{\max \left\{ \frac{m}{p}+\frac{1}{s}-1,0\right\} }.
\end{eqnarray*}

Now we prove the optimality. The optimality of the case $s\geq \frac{p}{p-m}$
is obvious. We just need to consider the case $s<\frac{p}{p-m}.$ Consider $%
A:\ell _{p}^{n}\times \cdots \times \ell _{p}^{n}\rightarrow \mathbb{K}$
(this idea we have borrowed from \cite{dimant}) given by%
\begin{equation*}
A(x^{(1)},...,x^{(m)})=\sum\limits_{j=1}^{n}x_{j}^{(1)}....x_{j}^{(m)}.
\end{equation*}%
Note that, from the H\"{o}lder inequality we obtain%
\begin{equation*}
\left\Vert A\right\Vert \leq n^{\frac{p-m}{p}}
\end{equation*}%
and thus, if the inequality (\ref{nt}) holds with $n^{t},$ we would have%
\begin{equation*}
n^{\frac{1}{s}}\leq Cn^{\frac{p-m}{p}}n^{t}
\end{equation*}%
and hence%
\begin{equation*}
t\geq \frac{m}{p}+\frac{1}{s}-1.
\end{equation*}

\textbf{Proof of (b). }For $p\geq 2,$ since $\ell _{p}$ has cotype $p$ with
cotype constant $1$ (the fact that the cotype constant is $1$ can be seen in
\cite[Lemma 2.3]{tonge} or \cite[page 29]{david}), we know from (\cite%
{botelho}) that every continuous $m$-linear form $T:\ell _{p}\times \cdots
\times \ell _{p}\rightarrow \mathbb{K}$ is absolutely $\left( \frac{p}{m}%
;1,....,1\right) $-summing and the summing norm of $T$ coincides with $%
\left\Vert T\right\Vert $. From the inclusion theorem for absolutely summing
multilinear operators (see \cite[Proposition 3.5]{matos} or \cite[%
Proposition 3.3]{david}) these forms are also absolutely $\left( 1;p^{\ast
},...,p^{\ast }\right) $-summing and, \textit{a fortiori}, absolutely $%
\left( s;p^{\ast },...,p^{\ast }\right) $-summing regardless of the $s\geq 1$%
. So we conclude that%
\begin{align*}
& \left( \sum\limits_{j=1}^{n}\left\vert T\left( e_{j},...,e_{j}\right)
\right\vert ^{s}\right) ^{\frac{1}{s}} \\
& \leq \left\Vert T\right\Vert \left( \sup_{\varphi \in B_{\left( \ell
_{p}^{n}\right) ^{\ast }}}\sum\limits_{j=1}^{n}\left\vert \varphi \left(
e_{j}\right) \right\vert ^{p^{\ast }}\right) ^{\frac{m}{p^{\ast }}} \\
& =\left\Vert T\right\Vert .
\end{align*}%
The optimality is immediate since we cannot have anything better than $%
n^{0}. $

\bigskip

\bigskip \textbf{Proof of (c). }For $1<p<2,$ since $\ell _{p}$ has cotype $2$
with cotype constant $1$ (again, the information that the cotype constant is
$1$ can be seen in \cite[Lemma 2.3]{tonge} or \cite[page 29]{david}), we
conclude from (\cite{botelho}) that every continuous $m$-linear form $T:\ell
_{p}\times \cdots \times \ell _{p}\rightarrow \mathbb{K}$ is absolutely $%
\left( \frac{2}{m};1,....,1\right) $-summing and, moreover, the respective
absolutely summing norm coincides with the $\sup $ norm. From the inclusion
theorem (see \cite[Proposition 3.5]{matos} or \cite[Prop. 3.3]{david}) for
absolutely summing multilinear operators, these forms are also absolutely $%
\left( s;\frac{2sm}{sm+2},...,\frac{2sm}{sm+2}\right) $-summing. Since $p<2$
it is easy to check that
\begin{equation}
\frac{2sm}{sm+2}<p^{\ast }.  \label{33}
\end{equation}%
Using (\ref{33}) we thus conclude that%
\begin{align*}
& \left( \sum\limits_{j=1}^{n}\left\vert T\left( e_{j},...,e_{j}\right)
\right\vert ^{s}\right) ^{\frac{1}{s}} \\
& \leq \left\Vert T\right\Vert \left[ \left( \sup_{\varphi \in B_{\left(
\ell _{p}^{n}\right) ^{\ast }}}\sum\limits_{j=1}^{n}\left\vert \varphi
\left( e_{j}\right) \right\vert ^{\frac{2sm}{sm+2}}\right) ^{\frac{sm+2}{2sm}%
}\right] ^{m} \\
& =\left\Vert T\right\Vert \left( \sum\limits_{j=1}^{n}\left( n^{\frac{-1}{%
p^{\ast }}}\right) ^{\frac{2sm}{sm+2}}\right) ^{\frac{sm+2}{2s}} \\
& =\left\Vert T\right\Vert \left( n\cdot \left( n^{\frac{1}{p}-1}\right) ^{%
\frac{2sm}{sm+2}}\right) ^{\frac{sm+2}{2s}} \\
& =\left\Vert T\right\Vert n^{\frac{sm+2}{2s}+\frac{m}{p}-m} \\
& =\left\Vert T\right\Vert n^{\frac{spm+2p+2ms-2spm}{2sp}} \\
& =\left\Vert T\right\Vert n^{\frac{2ms+2p-spm}{2sp}}.
\end{align*}

\bigskip The authors were kindly informed by Daniel Galicer and Pilar Rueda,
separately and in a private communication, and using different approaches, that the final solution to this
problem is the following:

\bigskip

\begin{theorem}
\label{ggg} Let $m,n$ be positive integers, $p_{1},...,p_{m}\geq 1$ and $s>0$%
. Let $C:=C(m,n,p_{1},...,p_{m},s)$ be the best constant that fulfils the
following inequality: for every $m$-linear form $A:\ell _{p_{1}}^{n}\times
\dots \times \ell _{p_{m}}^{n}\rightarrow \mathbb{K}$ we have
\begin{equation*}
\left( \sum_{j=1}^{n}|A(e_{j},\dots ,e_{j})|^{s}\right) ^{1/s}\leq C\;\Vert
A\Vert _{\mathcal{L}(\ell _{p_{1}}^{n},....,\ell _{p_{m}}^{n})}.
\end{equation*}%
The value of $C$ is exactly $%
\begin{cases}
(a)\;\;\;n^{1/s} & \text{ if }\frac{1}{p_{1}}+\cdots +\frac{1}{p_{m}}\geq 1,
\\
(b)\;\;n^{\max \left\{ \frac{1}{p_{1}}+\cdots +\frac{1}{p_{m}}+\frac{1}{s}%
-1,0\right\} } & \text{ if }\frac{1}{p_{1}}+\cdots +\frac{1}{p_{m}}\leq 1.%
\end{cases}%
$
\end{theorem}

\bigskip

\bigskip

\begin{remark}
As mentioned above, this final result is due to Daniel Galicer and Pilar Rueda
(independently, and using different approaches), so the present note will
not be submitted for publication.
\end{remark}

\end{document}